\documentstyle[12pt]{article}
\newcommand{\sect}[1]{\section{#1}\setcounter{equation}{0}}

\font\mbn=msbm10 scaled \magstep1
\font\mbs=msbm7 scaled \magstep1
\font\mbss=msbm5 scaled \magstep1
\newfam\mbff
\textfont\mbff=\mbn
\scriptfont\mbff=\mbs
\scriptscriptfont\mbff=\mbss\def\mbf{\fam\mbff}
\def\Re{{\mbf R}}

\def\Z{{\mbf Z}}
\def\Co{{\mbf C}}
\def\To{{\mbf T}}

\def\N{{\mbf N}}
\newtheorem{Th}{Theorem}[section]
\newtheorem{Lm}[Th]{Lemma}
\newtheorem{C}[Th]{Corollary}
\newtheorem{D}[Th]{Definition}
\newtheorem{Proposition}[Th]{Proposition}
\newtheorem{R}[Th]{Remark}
\newtheorem{E}[Th]{Example}
\author{Alexander Brudnyi\thanks{Research supported in part by NSERC.
}\\
Department of Mathematics and Statistics\\
University of Calgary, Calgary\\
Canada}
\title{Integral Representations of Holomorphic Functions on Coverings of
Pseudoconvex Domains in Stein Manifolds}
\date{} 
\begin{document} 
\maketitle
\begin{abstract}
{The classical integral representation formulas for holomorphic
functions defined on pseudoconvex domains in Stein manifolds play an 
important role in the 
constructive theory of functions of several complex variables.
In this paper we construct similar formulas for certain 
classes of holomorphic functions defined on coverings of such domains.}
\end{abstract}
%\tableofcontents
%===========================
\sect{\hspace*{-1em}. Introduction.}
{\bf 1.1.} The method of integral representations for holomorphic functions
works successfully in various problems of the theory of functions
of several complex variables, e.g., in problems of uniform estimates
for solutions of the Cauchy-Riemann equations, uniform estimates for
extensions of holomorphic functions from submanifolds, uniform approximation
of holomorphic functions that are continuous on the boundary, etc.
(We refer to the book of Henkin and Leiterer [HL] devoted to this
subject.)   
In the present paper we use this method to study holomorphic functions
of slow growth defined on unbranched coverings of pseudoconvex domains in
Stein manifolds. In particular, we will show that many known results
for holomorphic functions on such domains can be extended to 
similar results for holomorphic functions of slow 
growth defined on their coverings. Our  approach is based on new integral 
representation formulas for holomorphic functions of slow growth on 
coverings of pseudoconvex domains. These formulas generalize the classical 
Leray integral formula and certain of its developments. 
(Note that the classical 
integral formulas are usually applied to bounded domains with rectifiable 
boundaries, while infinite coverings of such domains may have not 
these properties.) The application of our integral formulas allows to 
reduce some problems for  holomorphic functions on the covering of a domain 
to analogous problems for Banach-valued holomorphic functions on the domain 
itself. In our proofs we exploit some
ideas previously used in [Br1], [Br2] in the area of the Corona problem and 
based on infinite-dimensional versions of Cartan's A and B theorems 
originally proved by Bungart [B].

In [Br3] we apply our technique to extend and strengthen
certain results of Gromov, Henkin and Shubin [GHS]
on holomorphic $L^{2}$-functions on coverings of pseudoconvex manifolds
in the case of coverings of Stein manifolds.

To formulate our main results we first introduce basic notation 
and definitions. Throughout this paper we consider complex manifolds $M$
satisfying the condition
\begin{equation}\label{e1}
M\subset\subset\widetilde M\subset N\ \ \ {\rm and}\ \ \ 
\pi_{1}(M)=\pi_{1}(N)
\end{equation}
where $M$ and $\widetilde M$ are open connected subsets of a complex 
manifold $N$, and $\widetilde M$ is Stein. 
(Here $\pi_{1}(X)$ denotes the fundamental group of $X$.) 

For instance, this condition is valid for $M$ a strictly pseudoconvex domain
or an analytic polyhedra in a Stein manifold.

Let $r:N'\to N$ be an unbranched covering of $N$. By $M'=r^{-1}(M)$ we denote
the corresponding unbranched covering of $M$. Condition (\ref{e1}) implies
that $M'$ is an open connected subset of $N'$ and 
$\pi_{1}(M')=\pi_{1}(N')$. Let
$\psi:N'\to\Re_{+}$ be such that $\log\psi$ is
uniformly continuous with respect to the path metric induced by a
Riemannian metric pulled back from $N$. We introduce the Banach space 
${\cal H}_{p,\psi}(M')$, $1\leq p\leq\infty$, of 
functions $f$ holomorphic on $M'$ with norm
\begin{equation}\label{e3}
|f|_{p,\psi}:=\sup_{x\in M}\left(\sum_{y\in r^{-1}(x)}
|f(y)|^{p}\psi(y)\right)^{1/p} .
\end{equation}
(Here ${\cal H}_{\infty,\psi}(M')$ is the Banach space of
bounded holomorphic functions on $M'$.) 
\begin{E}\label{ex1}
{\rm Let $d$ be the path metric on $N'$ obtained by
the pullback of a Riemannian metric defined on $N$. Fix
a point $o\in M'$ and set
$$
d_{o}(x):=d(o,x)\ ,\  \ \ x\in N'\ .
$$
It is easy to show by means of the triangle inequality that as the 
function $\psi$ one can take, e.g., $(1+d_{o})^{\alpha}$
or $e^{\alpha d_{o}}$ with $\alpha\in\Re$.
(For instance, if $N'$ is a strip \penalty-10000
$\{z=x+iy\ :\  |y|<L\}\subset\Co$ with
the action of group $\Z$ given by translations along $\Re$, i.e.,
$N'$ is a regular covering of an annulus, one can take
as $\psi$ either the functions $(1+|x|)^{\alpha}$ or $e^{\alpha |x|}$, 
$\alpha \in\Re$.)}
\end{E}
\begin{R}\label{r1}
{\rm Let $dV_{M'}$ be the Riemannian volume form on the covering $M'$ 
obtained by a Riemannian metric pulled back from $N$. 
Note that every $f\in {\cal H}_{p,\psi}(M')$ also belongs 
to the Banach space $H^{p}_{\psi}(M')$ of holomorphic functions $g$ on
$M'$ with norm
$$
\left(\int_{z\in M'}|g(z)|^{p}\psi(z)dV_{M'}(z)\right)
^{1/p} .
$$
Moreover, one has a continuous embedding 
${\cal H}_{p,\psi}(M')\hookrightarrow H^{p}_{\psi}(M')$.}
\end{R}

Let $x\in M$ and $r:M'\to M$ be an unbranched covering of $M$.
We introduce the
Banach space $l_{p,\psi, x}(M')$, $1\leq p\leq\infty$,
of functions $g$ on $r^{-1}(x)$ with norm
\begin{equation}\label{e5}
|g|_{p,\psi,x}:=\left(\sum_{y\in
r^{-1}(x)}|g(y)|^{p}\psi(y)\right)^{1/p} .
\end{equation}
Also, for Banach spaces $E$ and $F$, by 
${\cal B}(E,F)$ we denote the space of
all linear bounded operators $E\to F$ with norm $||\cdot||$.
\begin{Th}\label{te1}
Suppose that $M$ satisfies condition (\ref{e1}). Then for any
$p\in [1,\infty]$ there exists a family 
$\{L_{z}\in {\cal B}(l_{p,\psi,z}(M'),
{\cal H}_{p,\psi}(M'))\}_{z\in M}$ holomorphic 
in $z$ such that 
$$
(L_{z}h)(x)=h(x)\ \ \ {\rm for\ any}\ \ \ h\in l_{p,\psi,z}(M')\ \ \
{\rm and}\ \ \ x\in r^{-1}(z)\ .
$$
Moreover,
$$
\sup_{z\in M}||L_{z}||<\infty\ .
$$
\end{Th}
{\bf 1.2.} Theorem \ref{te1} allows to obtain integral representation
formulas for holomorphic functions from ${\cal H}_{p,\psi}(M')$ by 
means of known integral 
formulas for holomorphic functions on $M$. As an example we will show how to
get such formulas from the classical Leray integral formula,
the basis of many other integral formulas. We first recall this
formula itself.

For vectors $\xi ,\eta\in\Co^{n}$ we set
$$
<\eta,\xi>=\sum_{j=1}^{n}\eta_{j}\cdot\xi_{j}\ .
$$
Also, we set 
$$
\omega(\xi)=d\xi_{1}\wedge\cdots\wedge d\xi_{n}\ \ \ {\rm and}\ \ \
\omega'(\eta)=\sum_{k=1}^{n}(-1)^{k-1}\eta_{k}d\eta_{1}\wedge\cdots
\wedge d\eta_{k-1}\wedge d\eta_{k+1}\wedge\cdots\wedge d\eta_{n}\ .
$$
Let $M\subset\Co^{n}$ be a domain and $z\in M$ be a fixed point.
Consider in the domain $Q=\Co^{n}\times M$ with 
coordinates $\eta=(\eta_{1},\dots,\eta_{n})\in\Co^{n}$ and
$\xi=(\xi_{1},\dots,\xi_{n})\in M$ the hypersurface of the form
$$
P_{z}=\{(\eta,\xi)\in Q\ :\ <\eta,\ \xi-z>=0\}\ .
$$
Let $h_{z}$ be a $(2n-1)$-dimensional cycle in the domain $Q\setminus P_{z}$
such that its projection onto $M\setminus\{z\}$ is homologous to
$\partial M$. Then for any holomorphic function $f$ defined on $M$
we have (see [L])
\begin{equation}\label{le1}
f(z)=\frac{(n-1)!}{(2\pi i)^{n}}\int_{h_{z}}f(\xi)
\frac{\omega'(\eta)\wedge\omega(\xi)}{<\eta,\ \xi-z>^{n}}\ .
\end{equation}
Note that since the integral kernel in this formula is bounded and continuous
on $h_{z}$, similar formulas are valid for Banach-valued holomorphic functions
defined on $M$. 

Now from Theorem \ref{te1} we obtain 
\begin{C}\label{co1}
Suppose that $M\subset\Co^{n}$ satisfies condition (\ref{e1}). Then
under the assumptions of the Leray integral formula
for any holomorphic function $f\in {\cal H}_{p,\psi}(M')$ we have
\begin{equation}\label{le1'}
f(x)=\frac{(n-1)!}{(2\pi i)^{n}}\int_{h_{z}}L_{\xi}(f|_{r^{-1}(\xi)})(x)
\frac{\omega'(\eta)\wedge\omega(\xi)}{<\eta,\ \xi-z>^{n}}\ \ \ 
{\rm for\ any}\ \ \ x\in r^{-1}(z)\ .
\end{equation}
\end{C}

Similarly to (\ref{le1'}) we obtain extensions of other known 
integral formulas.
To get such an extension replace only the integrand function
$f(\xi)$ in an integral formula by $L_{\xi}(f|_{r^{-1}(\xi)})(x)$. In 
the same way one obtains multi-dimensional Cauchy-Green and
Koppelman-Leray type formulas for some classes of differential forms defined 
on coverings of $M$ satisfying (\ref{e1}). Several such formulas are
presented in the Appendix (see also [HL] for an exposition
of known integral formulas).\\  
{\bf 1.3.} Let us formulate some other applications of Theorem \ref{te1}.

Let $r: M'\to M$ be a covering.
Suppose that $\{x_{n}\}_{n\geq 1}\subset M$ converges to $x\in M$. Then
for sufficiently big $n$ we can arrange $r^{-1}(x_{n})$ and $r^{-1}(x)$ in 
sequences 
$\{y_{in}\}_{i\geq 1}$ and $\{y_{i}\}_{i\geq 1}$ such that every
$\{y_{in}\}$ converges to $y_{i}$ as $n\to\infty$. For such $n$ we
define maps $\tau_{n}(x): r^{-1}(x)\to r^{-1}(x_{n})$ so that 
$\tau_{n}(y_{i})=y_{in}$, $i\in\N$. Below, $\tau_{n}^{*}$ denotes
the transpose map generated by $\tau_{n}$ on functions defined on
$r^{-1}(x_{n})$ and $r^{-1}(x)$.
\begin{D}\label{cont1} 
Let $X\subset M$ be a subset. We say that a function $f$ on $r^{-1}(X)$
belongs to the class $C_{p,\psi}(r^{-1}(X))$ if
\begin{itemize}
\item[{\rm (1)}]
$f|_{r^{-1}(x)}\in l_{p,\psi,x}(M')$ for any $x\in X$ and 
\item[{\rm (2)}]
for any $x\in X$ and 
any sequence $\{x_{n}\}\subset X$ converging to $x$ the sequence of functions
$\{\tau_{n}^{*}(f|_{r^{-1}(x_{n})})\}$ converges to 
$f|_{r^{-1}(x)}$ in the norm of $l_{p,\psi,x}(M')$.
\end{itemize}
\end{D}

Suppose that $X$ belongs to a complex manifold $X_{a}\subset M$ with
$dim_{\Co}X_{a}\geq 1$ and any other complex manifold of
dimension $d$, $1\leq d<dim_{\Co}X_{a}$, does not contain $X$.
In what follows the interior of $r^{-1}(X)$ is defined with respect to
$r^{-1}(X_{a})$.
\begin{D}\label{holo}
By ${\cal H}_{p,\psi}(r^{-1}(X))$ we denote the Banach space 
of functions $f\in C_{p,\psi}(r^{-1}(X))$ holomorphic in interior 
points of $r^{-1}(X)$ with norm 
\begin{equation}\label{normx}
|f|_{p,\psi}^{X}:=\sup_{x\in X}|f|_{r^{-1}(x)}|_{p,\psi,x}\ .
\end{equation}
\end{D}
\begin{R}\label{bound}
{\rm (a) Arguments used in the proof of Proposition \ref{sec1} 
stated below show
that if $X\subset X_{a}$ is open then any 
function $f$ holomorphic on $r^{-1}(X)$ 
satisfying part (1) of Definition \ref{cont1} and such that 
$|f|_{p,\psi}^{X}<\infty$ belongs to ${\cal H}_{p,\psi}(r^{-1}(X))$.\\
(b) Suppose that $M$ satisfies condition (\ref{e1}) and that the function 
$\psi$ in the definition of ${\cal H}_{p,\psi}(M')$ is such that 
$1\in {\cal H}_{p,\psi}(M')$. Then it 
is easy to see that any function on $r^{-1}(X)$ 
uniformly continuous with respect to the path metric pulled back from
$N$ satisfies part (2) of Definition \ref{cont1}. In particular,
this is always true for $p=\infty$.}
\end{R}

Suppose that a domain $M\subset\Co^{N}$ satisfies (\ref{e1}). Let 
$D\subset\subset M$ be a 
strictly pseudoconvex open set (not necessarily $D\neq\emptyset$), and
$\rho$ be a strictly plurisubharmonic $C^{2}$-function in a neighbourhood
$O$ of $\partial D$ such that
\begin{equation}\label{comp}
D\cap O=\{z\in O\ : \rho(z)< 0\}\ \ \ {\rm and}\ \ \
N(\rho):=\{z\in O\ :\ \rho(z)=0\}\subset\subset O\ .
\end{equation}
We set $C:=D\cup N(\rho)$. 
\begin{Th}\label{te2}
There is a neighbourhood $\Omega$ of $C$ such that
every 
$f\in {\cal H}_{p,\psi}(r^{-1}(C))$ can be uniformly approximated in 
the norm of ${\cal H}_{p,\psi}(r^{-1}(C))$ by functions from
${\cal H}_{p,\psi}(r^{-1}(\Omega))$.
\end{Th}
\begin{E}\label{real}
{\rm Let $C$ be a (connected) compact real-analytic manifold.
By the Grauert theorem [Gr] we can assume that $C$ is an analytic 
submanifold of some $\Re^{n}$. 
It is easy to show that $C$ is the zero set of a non-negative strictly 
plurisubharmonic $C^{2}$-function defined in a neighbourhood of
$C$. Let $C_{\epsilon}$ be an $\epsilon$-neighbourhood of $C$ in 
$\Co^{n}$ where $\epsilon$ is sufficiently small. 
Then $C_{\epsilon}$ satisfies condition (\ref{e1}).
Let $r:C'\to C$ be a covering of $C$, and 
$r:C_{\epsilon}'\to C_{\epsilon}$ be the 
covering of $C_{\epsilon}$ such that $C'\subset C_{\epsilon}'$ 
(i.e., for this covering $\pi_{1}(C')=\pi_{1}(C_{\epsilon}')$).  
Now, Theorem \ref{te2} implies that for a sufficiently small $\epsilon$ 
any $C_{p,\psi}$-function on $C'$ can be uniformly approximated
in the norm of ${\cal H}_{p,\psi}(C')$ by functions from
${\cal H}_{p,\psi}(C_{\epsilon}')$.}
\end{E} 

Finally, let us formulate a result on bounded extension of holomorphic 
functions from complex submanifolds.

Suppose that $M$ is a strictly pseudoconvex open set (with not necessarily
smooth boundary) in a Stein manifold satisfying
condition (\ref{e1}). Let $Y$ be a closed complex 
submanifold of some neighbourhood of $\overline{M}$. Consider a covering 
$r:M'\to M$. Since $M$ satisfies 
condition (\ref{e1}), there exists a covering 
$r:N'\to N$ such that $M'$ is an open subset of
$N'$. By $\overline{M'}$ we denote the closure of $M'$ in $N'$.
Also, we set $X:=Y\cap M$.

\begin{Th}\label{te3}
\begin{itemize}
\item[(1)]
For every function $f\in {\cal H}_{p,\psi}(r^{-1}(X))$, there exists
a function $F\in {\cal H}_{p,\psi}(M')$ such that $F=f$ on $r^{-1}(X)$.
\item[(2)]
For every function $f\in {\cal H}_{p,\psi}(r^{-1}(\overline{X}))$, 
there exists a function
$F\in {\cal H}_{p,\psi}(\overline{M'})$ such that $F=f$ on 
$r^{-1}(\overline{X})$.
\end{itemize}
\end{Th}
(Recall that ${\cal H}_{p,\psi}(r^{-1}(\overline{X}))$
is the space of $C_{p,\psi}$-functions on $r^{-1}(\overline{X})$
holomorphic in $r^{-1}(X)$ with the norm defined by (\ref{normx}).)

Analogous results hold for a connected component $Z$
of $r^{-1}(X)$. In this case we define ${\cal H}_{p,\psi}(Z)$
as the space of functions on $Z$ whose extensions to 
$r^{-1}(X)$ by 0 belong to ${\cal H}_{p,\psi}(r^{-1}(X))$. One
defines ${\cal H}_{p,\psi}(\overline{Z})$ similarly.

%=========================
\sect{\hspace*{-1em}. Auxiliary Results.}
For the standard facts about bundles see, e.g., Hirzebruch's book [Hi]. 
In what follows, all topological spaces are assumed to be finite or 
infinite dimensional.\\
{\bf 2.1.} Let $X$ be a complex analytic space and $S$ be a complex 
analytic Lie group with
unit $e\in S$. Consider an effective holomorphic action of $S$ on
a complex analytic space $F$. Here {\em holomorphic action} means a 
holomorphic map $S\times F\rightarrow F$ sending 
$s\times f\in S\times F$ to
$sf\in F$ such that $s_{1}(s_{2}f)=(s_{1}s_{2})f$ and $ef=f$ for any
$f\in F$. {\em Efficiency} means that the condition $sf=f$ for some $s$ and 
any $f$ implies that $s=e$. 
\begin{D}\label{de2}
A complex analytic space $W$ together with a holomorphic map (projection)
$\pi:W\rightarrow X$ is called a holomorphic bundle on $X$ with 
structure group $S$ and fibre $F$, if there exists a
system of coordinate transformations, i.e., if
\begin{itemize}
\item[{\rm (1)}]
there is an open cover ${\cal U}=\{U_{i}\}_{i\in I}$ of $X$ and a
family of biholomorphisms 
$h_{i}:\pi^{-1}(U_{i})\rightarrow U_{i}\times F$,
that map ``fibres'' $\pi^{-1}(u)$ onto $u\times F$;
\item[{\rm (2)}] 
for any $i,j\in I$ there are elements 
$s_{ij}\in {\cal O}(U_{i}\cap U_{j}, S)$ such that
$$
(h_{i}h_{j}^{-1})(u\times f)=u\times s_{ij}(u)f\ \ \ {\rm for\ any}\ \
u\in U_{i}\cap U_{j},\ f\in F\ .
$$
\end{itemize}
In particular, a 
holomorphic bundle $\pi:W\rightarrow X$ whose fibre is a Banach space 
$F$ and the structure group is
$GL(F)$ (the group of linear invertible transformations of $F$) is
called a holomorphic Banach vector bundle. 

A holomorphic section of a holomorphic bundle $\pi:W\rightarrow X$
is a holomorphic map $s:X\rightarrow W$ satisfying $\pi\circ s=id$.

Let $\pi_{i}:W_{i}\rightarrow X$, $i=1,2$, be holomorphic Banach vector
bundles. A holomorphic map $f:W_{1}\rightarrow W_{2}$ satisfying
\begin{itemize}
\item[{\rm (a)}]
$f(\pi_{1}^{-1}(x))\subset\pi_{2}^{-1}(x)$ for any $x\in X$;
\item[{\rm (b)}] 
$f|_{\pi_{1}^{-1}(x)}$ is a linear continuous map of the corresponding
Banach spaces,
\end{itemize}
is called a homomorphism. If, in addition, $f$ is a 
homeomorphism, then $f$ is called an isomorphism.
\end{D}
We will use
the following construction of holomorphic bundles 
(see, e.g. [Hi,\ Ch.1]):

Let $S$ be a complex analytic Lie group and 
${\cal U}=\{U_{i}\}_{i\in I}$ be an open cover of $X$. 
By $Z_{\cal O}^{1}({\cal U},S)$ we denote the set of holomorphic $S$-valued 
${\cal U}$-cocycles. By definition, 
$s=\{s_{ij}\}\in Z_{\cal O}^{1}({\cal U},S)$, 
where $s_{ij}\in {\cal O}(U_{i}\cap U_{j}, S)$ and
$s_{ij}s_{jk}=s_{ik}\ \ \  {\rm on}\ \ \  U_{i}\cap U_{j}\cap U_{k}$.
Consider the disjoint union $\sqcup_{i\in I}U_{i}\times F$ and for any 
$u\in U_{i}\cap U_{j}$ identify the point $u\times f\in U_{j}\times F$ with 
$u\times s_{ij}(u)f\in U_{i}\times F$.
We obtain a holomorphic bundle $W_{s}$ on $X$ whose projection is induced 
by the projection $U_{i}\times F\rightarrow U_{i}$.
Moreover, any holomorphic bundle on $X$ with structure group $S$ and 
fibre $F$ is isomorphic (in the category of holomorphic bundles) to a 
bundle $W_{s}$. 
\begin{E}\label{constbun}
{\rm {\bf (a)}  Let $M$ be a complex manifold.
For any subgroup $H\subset\pi_{1}(M)$ consider the unbranched covering
$r:M(H)\rightarrow M$ corresponding to $H$. We will describe $M(H)$
as a holomorphic bundle on $M$.

First, assume that $H\subset\pi_{1}(M)$ is a normal subgroup.
Then $M(H)$ is a regular covering of $M$ and the quotient group
$G:=\pi_{1}(M)/H$ acts holomorphically on $M(H)$ by deck 
transformations.  It is well known that $M(H)$ in this case 
can be thought of as a {\em principle fibre bundle} on $M$ with fibre
$G$ (here $G$ is equipped with the discrete topology). 
Namely, let us consider the map $R_{G}(g):G\rightarrow G$, $g\in G$,
defined by the formula
$$
R_{G}(g)(q)=q\cdot g^{-1},\ \ \ q\in G.
$$ 
Then there is an open cover ${\cal U}=\{U_{i}\}_{i\in I}$ of $M$ by sets
biholomorphic to open Euclidean balls in some $\Co^{n}$ and a locally constant
cocycle $c=\{c_{ij}\}\in Z_{\cal O}^{1}({\cal U},G)$ 
such that $M(H)$ is biholomorphic
to the quotient space of the disjoint union 
$V=\sqcup_{i\in I}U_{i}\times G$ by the equivalence relation:
$U_{i}\times G\ni x\times R_{G}(c_{ij})(q)\sim x\times q\in 
U_{j}\times G$.
The identification space is a holomorphic bundle with projection 
$r:M(H)\rightarrow M$ induced by the
projections $U_{i}\times G\rightarrow U_{i}$.
In particular, when $H=e$ we obtain the definition of the universal
covering $M_{u}$ of $M$. 

Assume now that $H\subset\pi_{1}(M)$ is not necessarily normal.
Let $X_{H}=\pi_{1}(M)/H$ be the set of cosets with respect
to the (left) action of $H$ on $\pi_{1}(M)$ defined by left multiplications. 
By $[Hq]\in X_{H}$ we denote the coset containing $q\in\pi_{1}(M)$.
Let $A(X_{H})$ be the group of all homeomorphisms of $X_{H}$
(equipped with the discrete topology). We define the homomorphism
$\tau:\pi_{1}(M)\rightarrow A(X_{H})$ by the formula:
$$
\tau(g)([Hq]):=[Hqg^{-1}],\ \ \ q\in\pi_{1}(M).
$$
Set $Q(H):=\pi_{1}(M)/Ker(\tau)$ and let 
$\widetilde g$ be the image of $g\in\pi_{1}(M)$ in $Q(H)$. 
By $\tau_{Q(H)}:Q(H)\rightarrow A(X_{H})$ we denote the unique 
homomorphism whose pullback to $\pi_{1}(M)$ coincides with $\tau$.
Consider the action of
$H$ on $V=\sqcup_{i\in I}U_{i}\times \pi_{1}(M)$ induced by the left 
action of $H$ on $\pi_{1}(M)$ and let
$V_{H}=\sqcup_{i\in I}U_{i}\times X_{H}$ be the corresponding quotient set.
Define the equivalence relation
$U_{i}\times X_{H}\ni x\times \tau_{Q(H)}(\widetilde c_{ij})(h)
\sim x\times h\in U_{j}\times X_{H}$ 
with the same $\{c_{ij}\}$ as in the definition of $M(e)$.
The corresponding quotient space is a holomorphic bundle with fibre $X_{H}$ 
biholomorphic to $M(H)$. \\
{\bf (b)}\ We retain the notation of example (a). Let 
$B$ be a complex Banach space and
$GL(B)$ be the group of invertible
bounded linear operators $B\to B$. Consider a homomorphism 
$\rho: G\rightarrow GL(B)$. Without loss of generality we assume that 
$Ker(\rho)=e$, for otherwise we can pass to the corresponding quotient
group. The {\em holomorphic Banach vector bundle
$E_{\rho}\rightarrow M$ associated with $\rho$ } is defined as 
the quotient of $\sqcup_{i\in I} U_{i}\times B$ by the equivalence
relation
$U_{i}\times B\ni x\times\rho(c_{ij})(w)\sim x\times w\in U_{j}\times B$
for any $x\in U_{i}\cap U_{j}$. Let us illustrate this construction by 
an example. 

Let $X_{H}=\pi_{1}(M)/H$ be the fibre of the covering $r:M(H)\to M$.
Consider a function $\phi:X_{H}\to\Re_{+}$ satisfying
\begin{equation}\label{dilat}
\phi(\tau_{Q(H)}(h)(x))\leq c_{h}\phi(x)\ ,\ \ \ h\in Q(H)\ ,\ x\in X_{H}\ ,
\end{equation}
where $c_{h}$ is a constant depending on $h$. By $l_{p,\phi}(X_{H})$, 
$1\leq p\leq\infty$, we denote
the Banach space of complex functions $f$ on $X_{H}$ with norm
\begin{equation}\label{p1}
||f||_{p,\phi}:=\left(\sum_{x\in X_{H}}|f(x)|^{p}\phi(x)\right)^{1/p} .
\end{equation}
Then according to (\ref{dilat}) the map $\rho$ defined by the formula
$$
[\rho(h)(f)](x):=f(\tau_{Q(H)}(h)(x))\ ,\ \ \ h\in Q(H)\ ,\ x\in X_{H}\ ,
$$ 
is a homomorphism of
$Q(H)$ into $GL(l_{p,\phi}(X_{H}))$. By $E_{p,\phi}(X_{H})$ we denote the 
holomorphic Banach vector bundle associated with this $\rho$.}
\end{E}
{\bf 2.2.} Let $r:M'\to M$ be a covering where $M'=M(H)$ (i.e.,
$\pi_{1}(M')=H$). In this part we establish a connection between Banach
spaces ${\cal H}_{p,\psi}(M')$ defined in section 1.1 and certain
spaces of holomorphic sections of bundles $E_{p,\phi}(X_{H})$.

We retain the notation of Example \ref{constbun}. 
Assume that $M$ satisfies condition (\ref{e1}), i.e., 
$M\subset\subset N$ and $\pi_{1}(M)=\pi_{1}(N)$.
Then $M(H)$ is embedded into $N(H)$.
(Without loss of generality we consider $M(H)$ as an
open subset of $N(H)$.) Let $\{V_{i}\}_{i\in I}$ be
a finite acyclic open cover of $\overline{M}$ by relatively compact sets.
We set $U_{i}:=V_{i}\cap M$ and consider the open cover 
${\cal U}=\{U_{i}\}_{i\in I}$ of $M$. Then as in Example \ref{constbun} (a) 
we can define $M(H)$ by a cocycle
$\widetilde c=\{\widetilde c_{ij}\}\in Z^{1}_{{\cal O}}({\cal U},Q(H))$
where $Q(H)=\pi_{1}(M)/Ker\ \!\tau$.

Further, let $\psi:N(H)\to\Re_{+}$ be a function such that $\log\psi$ is 
uniformly continuous with respect to the path metric induced by a Riemannian 
metric pulled back from $N$. Fix a point $z_{0}\in M$. Identifying
$r^{-1}(z_{0})$ with $X_{H}$ define the function $\phi: X_{H}\to\Re_{+}$ by 
the formula
$$
\phi(x):=\psi(x)\ ,\ \ \ x\in X_{H}\ (=r^{-1}(z_{0}))\ .
$$
\begin{Lm}\label{lem1}
The function $\phi$ satisfies inequality (\ref{dilat}).
\end{Lm}
{\bf Proof.}
We recall some facts from the theory of covering spaces (see, e.g., 
[Hu, Chapter III]. Let $h\in Q(H)$. Then there exists a closed path 
$\gamma\subset M$ passing through the point $z_{0}$ such that its 
lifting $\gamma_{x}\subset M(H)$ with the initial point $x\in r^{-1}(z_{0})$
has the endpoint $\tau_{Q(H)}(h)(x)$. By the definition of the metric on
$M(H)$, the length of every such $\gamma_{x}$ is the same.
This and uniform continuity of $\log\psi$ with respect to the path
metric defined by a Riemannian metric pulled back from $N$ imply that
every $\gamma_{x}$ can be covered by $k$ metric balls
$\{N_{j}(x)\}_{j=1}^{k}$ such that
$$
\frac{1}{2}\psi(y)\leq\psi(z)\leq 2\psi(y)\ \ \ {\rm for\ any}\ \ \
y,z\in N_{j}(x), \ \ \ 1\leq j\leq k\ .
$$
In particular, we have
$$
\frac{\psi(\tau_{Q(H)}(h)(x))}{\psi(x)}\leq 2^{k}\ \ \ {\rm for \ every}\ 
x\in r^{-1}(z_{0})\ . \ \ \ \ \ \Box
$$

According to this lemma, the bundle $E_{p,\phi}(X_{H})$ is well defined.
By definition, any holomorphic section of this bundle can be determined 
by the family $\{f_{i}(z,x)\}_{i\in I}$ of holomorphic functions on
$U_{i}$ with values in $l_{p,\phi}(X_{H})$ satisfying
$$
f_{i}(z,\tau_{Q(H)}(\widetilde c_{ij})(x))=f_{j}(z, x)\ \ \ {\rm for\ any}
\ \ \ z\in U_{i}\cap U_{j}\ .
$$
We introduce the Banach space $B_{p,\phi}(X_{H})$ of {\em
bounded holomorphic sections}
$f=\{f_{i}\}_{i\in I}$ of $E_{p,\phi}(X_{H})$ with norm
\begin{equation}\label{se1}
|f|_{p,\phi}:=\sup_{i\in I, z\in U_{i}}||f_{i}(z,\cdot)||_{p,\phi}\ .
\end{equation}
(Here $||\cdot ||_{p,\phi}$ is the norm on $l_{p,\phi}(X_{H})$, 
see (\ref{p1}).)

Further, let $f\in {\cal H}_{p,\psi}(M(H))$ (see section 1.1 for the 
definition). We identify $M(H)$ with the quotient set of
$V_{H}=\sqcup_{i\in I}U_{i}\times X_{H}$ as in Example \ref{constbun} (a).
This gives local coordinates on $M(H)$. Using these coordinates
we define the family
$\{f_{i}\}_{i\in I}$ of functions on $U_{i}$ with values in the space of 
functions on $X_{H}$ by the formula
\begin{equation}\label{cor1}
f_{i}(z,x):=f(z,x)\ ,\ \ \ z\in U_{i}\ ,\ i\in I\ ,\ x\in X_{H}\ .
\end{equation}
Then the following result holds.
\begin{Proposition}\label{sec1}
The correspondence $f\mapsto \{f_{i}\}_{i\in I}$ determines an isomorphism
of Banach spaces $D:{\cal H}_{p,\psi}(M(H))\to B_{p,\phi}(X_{H})$.
\end{Proposition}
{\bf Proof.} 
First we will prove that every $f_{i}(z,x)=f(z,x)$ determines
a holomorphic map of $U_{i}$ into $l_{p, \phi}(X_{H})$. Since by definition
$f_{i}(z,\tau_{Q(H)}(\widetilde c_{ij}(x))=f_{j}(z,x)$ for any 
$z\in U_{i}\cap U_{j}$, this will
show that $\{f_{i}\}$ is a holomorphic section of $E_{p,\phi}(X_{H})$.

Check that $f_{i}(z,\cdot))\in l_{p,\phi}(X_{H})$ for any $z\in U_{i}$. 
We will assume that $z_{0}$ in the definition of $\phi$ from Lemma
\ref{lem1} belongs to some $U_{k}$ and so $r^{-1}(z_{0})=z_{0}\times X_{H}$
in the above coordinates on $r^{-1}(U_{k})\subset M(H)$. Let
us fix some points $z_{i}\in U_{i}$. Since $M$ is compactly embedded into 
$N$ and the path metric on $N(H)$ is
defined by a Riemannian metric pulled back from $N$,
for every $(z,x)\in r^{-1}(U_{i})$, $i\in I$, there exists a path 
$\gamma_{i}(x)$ joining 
$(z_{i},x)$ with $(z,x)$ such that the lengths of all these paths are bounded
from above by a constant $L$. Further, by the covering homotopy theorem
for every $(z_{0},x)\in r^{-1}(z_{0})$ there exist a point 
$(z_{i},s_{i}(x))\in r^{-1}(U_{i})$ and a path $\widetilde\gamma_{i}(x)$
joining these two points in $M(H)$ such that 
$s_{i}(x)=\tau_{Q(H)}(h_{i})(x)$ for some $h_{i}\in Q(H)$ 
(independent of $x$) and the lengths of all $\widetilde\gamma_{i}(x)$
(for all $i\in I$ and $x\in X_{H}$) are bounded from above by $L_{1}$.
Since $\log\psi$ is uniformly continuous with respect to the path metric
obtained by a Riemannian metric pulled back from $N$, from here arguing
as in the proof of Lemma \ref{lem1} and using the boundedness of lengths
of all the paths $\widetilde\gamma_{i}(x)$ and $\gamma_{i}(x)$ we obtain 
that there is a constant $C$ depending on $M$ and $z_{0}$ such that
\begin{equation}\label{comp1}
\frac{1}{C}\psi(z,x)\leq\phi(x)\leq C\psi(z,x)\ \ \ 
{\rm for\ any}\ \ \ z\in U_{i},\ i\in I,\ x\in X_{H}\ .
\end{equation}
Now by the definition of the norm on ${\cal H}_{p,\psi}(M(H))$ we have
$$
\sup_{i\in I,z\in U_{i}}
\left(\sum_{x\in X_{H}}|f(z,x)|^{p}\psi(z,x)\right)^{1/p}=
|f|_{p,\psi}\ .
$$
Combining this with (\ref{comp1}) we get
\begin{equation}\label{comp2}
(1/C)^{1/p}\cdot|f|_{p,\psi}\leq
\sup_{i\in I,z\in U_{i}}||f(z,\cdot)||_{p,\phi}
\leq C^{1/p}\cdot |f|_{p,\psi}\ .
\end{equation}
This shows that every $f_{i}(z,\cdot)\in l_{p,\phi}(X_{H})$. 

Let us prove now that $f_{i}(\cdot,\cdot):U_{i}\to l_{p,\phi}(X_{H})$ is 
holomorphic. Without loss of generality we identify $U_{i}$ with an open 
subset of a certain $\Co^{n}$. Let $P\subset\subset U_{i}$,
$P=\{(z_{1},\dots,z_{n})\in\Co^{n}\ :\ 
|z_{j}-x_{j}|\leq r,\ 1\leq j\leq n\}$, be a polydisk with the center
at a point $x=(x_{1},\dots,x_{n})\in U_{i}$ and of radius $r$. 
By $P_{t}$ we denote its boundary torus (i.e., $P_{t}=\{z\in P\ :\
|z_{j}-x_{j}|=r,\ 1\leq j\leq n\}$). Let us introduce
$c_{\alpha_{1},\dots,\alpha_{n}}(y)$ by the formula
$$
c_{\alpha_{1},\dots,\alpha_{n}}(y)=\left(\frac{1}{2\pi i}\right)^{n}
\int_{P_{t}}
\frac{f_{i}(w,y)}{w_{1}^{\alpha_{1}+1}\cdots w_{n}^{\alpha_{n}+1}}
dw_{1}\cdots dw_{n}\ .
$$
We will show that 
\begin{equation}\label{comp3}
c_{\alpha_{1},\dots,\alpha_{n}}\in l_{p,\phi}(X_{H}) \ \ \ {\rm and}
\ \ \
||c_{\alpha_{1},\dots,\alpha_{n}}||_{p,\phi}\leq 
\frac{C^{1/p}|f|_{p,\psi}}{r^{\alpha_{1}+\dots +\alpha_{n}}}\ .
\end{equation}
In what follows $\To^{n}$ denotes the standard torus in $\Re^{n}$ with 
coordinates $t=(t_{1},\dots, t_{n})$.
Now, using the definition of $c_{\alpha_{1},\dots,\alpha_{n}}(y)$,
the triangle inequality for the norm $||\cdot||_{p,\phi}$, inequality
(\ref{comp2}) and the Lebesgue monotone convergence theorem we have
$$
\begin{array}{c}
\displaystyle
||c_{\alpha_{1},\dots,\alpha_{n}}||_{p,\phi}:=
\left(\sum_{y\in X_{H}}
|c_{\alpha_{1},\dots,\alpha_{n}}(y)|^{p}\phi(y)\right)^{1/p}\leq\\
\\
\displaystyle
\left(\sum_{y\in X_{H}}\left(\left(\frac{1}{2\pi}\right)^{n}
\int_{\To^{n}}\frac{|f_{i}(x+rt,y)|}{r^{\alpha_{1}+\dots +\alpha_{n}}}
\ \!dt_{1}\cdots dt_{n}\right)^{p}\cdot\phi(y)\right)^{1/p}\leq\\
\\
\displaystyle
\frac{1}{r^{\alpha_{1}+\dots +\alpha_{n}}}\cdot\left(\frac{1}{2\pi}\right)^{n}
\int_{\To^{n}}\left(\sum_{y\in X_{H}}|f_{i}(x+rt,y)|^{p}\phi(y)\right)^{1/p}
dt_{1}\cdots dt_{n}\leq 
\frac{C^{1/p}|f|_{p,\psi}}{r^{\alpha_{1}+\dots +\alpha_{n}}}\ .
\end{array}
$$
Inequality (\ref{comp3}) implies, in particular, that
$$
\sum_{\alpha_{1}+\dots+\alpha_{n}=0}^{\infty}
c_{\alpha_{1},\dots,\alpha_{n}}(z_{1}-x_{1})^{\alpha_{1}}\cdots
(z_{n}-x_{n})^{\alpha_{n}}
$$
converges uniformly in any polydisk $P'$ centered at $x$ of radius $<r$ to
a holomorphic $l_{p,\phi}(X_{H})$-valued function $F$. By definition, 
evaluation
of this function at every point $y\in X_{H}$ coincides with 
$f_{i}(\cdot, y)$
on any such $P'$. Thus $F(z)=f_{i}(z,\cdot)$ for $z\in P'$. Since
$x$ is an arbitrary point of $U_{i}$, this shows
that $f_{i}(\cdot,\cdot): U_{i}\to l_{p,\phi}(X_{H})$ is holomorphic.

Therefore by (\ref{comp2}) we obtain that the map 
$D:{\cal H}_{p,\psi}(M(H))\to B_{p,\phi}(X_{H})$, \penalty-10000
$f\mapsto\{f_{i}\}_{i\in I}$, is
well defined and continuous.

Conversely, let $\{f_{i}(z,x)\}_{i\in I}$ be a family of holomorphic 
$l_{p,\phi}(X_{H})$-valued functions on the cover $\{U_{i}\}_{i\in I}$ 
determining an element from $B_{p,\phi}(X_{H})$. We set
$$
f(z,x):=f_{i}(z,x)\ ,\ \ \ z\in U_{i}\ ,\ x\in X_{H}\ .
$$
Then by the definition $f$ is a holomorphic function on $M(H)$.
Now inequality (\ref{comp2}) easily implies that 
$|f|_{p,\psi}\leq C^{1/p}\cdot |f|_{p,\phi}$.\ \ \ \ \ $\Box$\\

The map $D$ from Proposition \ref{sec1} will be called the {\em direct
image map.} 
\begin{R}\label{re3}
{\rm One can easily see that the direct image map $D$ is 
an isometry in the case $\psi\equiv 1$. }
\end{R}
%=============
\sect{\hspace*{-1em}. Proof of Theorem \ref{te1}.}
{\bf 3.1.} Assume that $M$ satisfies condition (\ref{e1}).
Let $r:M'\to M$ be an unbranched covering.
By $R_{x}$ we denote the restriction map of functions from 
${\cal H}_{p,\psi}(M')$ to 
the fibre $r^{-1}(x)$, $x\in M$. The definitions of the norms on
${\cal H}_{p,\psi}(M')$ and $l_{p,\psi,x}(M')$ (see section 1.1)
imply that $R_{x}: {\cal H}_{p,\psi}(M')\to l_{p,\psi,x}(M')$ is a 
linear continuous operator.
\begin{Proposition}\label{pr1}
For every $x\in M$ there exists a linear continuous operator \penalty-10000
$C_{x}: l_{p,\psi,x}(M')\to {\cal H}_{p,\psi}(M')$ such that 
$$
R_{x}\circ C_{x}= id\ .
$$
\end{Proposition}
{\bf Proof.} As before we think of $M'$ as an open subset of 
$N'$ where $r:N'\to N$ is a covering such that $\pi_{1}(N')=\pi_{1}(M')$.
We also set $\widetilde M':=r^{-1}(\widetilde M)$. (Here by our 
assumption $M\subset\subset\widetilde M\subset N$ and $\widetilde M$ is 
Stein.) Recall that the fibre of the covering $N'$ is
$X_{H}:=\pi_{1}(N)/H$ where $H:=\pi_{1}(N')$ and the quotient is taken with 
respect to the
left action of $H$ on $\pi_{1}(N)$. In what follows,
according to Proposition \ref{sec1},
we identify ${\cal H}_{p,\psi}(M')$ with the space $B_{p,\phi}(X_{H})$
of bounded holomorphic sections of the Banach vector bundle
$E_{p,\phi}(X_{H})\to M$ where $\phi=\psi|_{r^{-1}(z_{0})}$ for some
$z_{0}\in M$ (see Example \ref{constbun}). 

Let us introduce the Banach space $B$ of complex functions $F$ 
defined on $X_{H}\times X_{H}$ with norm
$$
|F(h,g)|_{B}:=\max\left\{\sup_{h\in X_{H}}\left(
\sum_{g\in X_{H}}|F(h,g)|\cdot\frac{\phi(g)}{\phi(h)}\right)\ ,\ 
\sup_{g\in X_{H}}\left(\sum_{h\in X_{H}}|F(h,g)|
\right)\right\}
\ .
$$
It is easy to see by means of (\ref{dilat}) that the formula
$$
[\rho(q,s)(F)](h,g):=F(\tau_{Q(H)}(q)(h),\tau_{Q(H)}(s)(g)),\ \ \
h,g\in X_{H}\ ,\ q, s \in Q(H)\ ,
$$
determines a homomorphism of the group $Q(H)\times Q(H)$ into 
$GL(B)$ (see Example \ref{constbun} (a) for the definitions
of $\tau_{Q(H)}$ and $Q(H)$). Note that $Q(H)\times Q(H)$ is the quotient of 
$\pi_{1}(N\times N)=\pi_{1}(N)\times\pi_{1}(N)$. Thus the
associated with $\rho$ holomorphic Banach vector bundle 
$t:E_{\rho}\to N\times N$ on $N\times N$ is defined.
We identify the fibre $t^{-1}(x\times x)$ with $B$.
Further, let $\delta_{h}$ be a function on $X_{H}$ such that
$\delta_{h}(g)=1$ if $h=g$ and $\delta_{h}(g)=0$ if $h\neq g$.
We define a function $\Delta$ on $X_{H}\times X_{H}$ by the formula
$$
\Delta(h,g):=\delta_{h}(g)\ .
$$
Check that $\Delta\in B\ (=t^{-1}(x\times x))$. Indeed,
$$
|\Delta|_{B}:=\max\left\{\sup_{h\in X_{H}}\left(
\sum_{g\in X_{H}}|\delta_{h}(g)|\cdot\frac{\phi(g)}{\phi(h)}\ \right) ,
\ \sup_{g\in X_{H}}
\left(\sum_{h\in X_{H}}|\delta_{h}(g)|\right)\right\}=1\ .
$$
Now, according to Bungart [B, Lemma 3.3]) there is a
holomorphic section $S$ of $E_{\rho}$ such that $S(x\times x)=\Delta$.
Since $\overline{M}\subset N$ is compact,  
$S|_{M\times M}$ is a bounded holomorphic section of 
$E_{\rho}|_{M\times M}$, cf. (\ref{se1}).

Let $\{V_{i}\}_{i\in I}$ be a finite acyclic open cover of 
$\overline{M}$. We
set $U_{i}=V_{i}\cap M$ and consider the open cover 
${\cal U}=\{U_{i}\times U_{j}\}_{i,j\in I}$ of $M\times M$. 
According to Example \ref{constbun} (a) $M'$ is defined on the cover 
$\{U_{i}\}_{i\in I}$ of $M$ by a cocycle 
$\widetilde c=\{\widetilde c_{ij}\}\in 
{\cal Z}^{1}_{{\cal O}}(\{U_{i}\},Q(H))$.
Then $M'\times M'$ is defined on ${\cal U}$ by the cocycle
$\widetilde c\times\widetilde c=\{\widetilde c_{ij}\times\widetilde c_{ij}\}
\in {\cal Z}^{1}_{{\cal O}}({\cal U},Q(H)\times Q(H))$. Using this
construction one represents the restriction of the section $S$ to 
${\cal U}$ by a family $\{S_{ij}(z,h,w,g)\ :\ h,g\in X_{H}\ ,\ 
z\times w\in U_{i}\times U_{j}\}_{i,j\in I}$ 
of holomorphic functions on $U_{i}\times U_{j}$ with values in $B$ satisfying
for $z\times w\in (U_{i}\times U_{j})\cap (U_{l}\times U_{m})$
\begin{equation}\label{intersect}
S_{ij}(z,\tau_{Q(H)}(\widetilde c_{il})(h),w,
\tau_{Q(H)}(\widetilde c_{mj})(g))=S_{lm}(z,h,w,g)\ .
\end{equation}
Since $S$ is bounded on $M\times M$,
\begin{equation}\label{bou1}
\sup_{i,j\in I,z\times w\in U_{i}\times U_{j}}|S_{ij}(z,\cdot,w,\cdot)|_{B}=
C<\infty\ .
\end{equation}
Also, by the definition we have for $x\times x\in M\times M$ (say, e.g.,
$x\in U_{i_{0}}$, $i_{0}\in I$) 
\begin{equation}\label{rest1}
S_{i_{0}i_{0}}(x,h,x,g)=\delta_{h}(g)\ .
\end{equation}
Next, consider the restriction of $S$ to $x\times M$ and the open
cover  $\{U_{i_{0}}\times U_{i}\}_{i\in I}$ of $x\times M$.
We define a family $\{T_{i}\}_{i\in I}$ of linear operators
on the set of complex functions $v$ on $X_{H}$ by the formula
\begin{equation}\label{oper}
T_{i}(v)(z,g):=\sum_{h\in X_{H}}v(h)S_{i_{0}i}(x,h,z,g)\ ,\ \ \ z\in U_{i},\ 
g\in X_{H}\ .
\end{equation}
\begin{Lm}\label{bounded}
(a)
$T_{i}$ is a bounded operator from 
$l_{p,\phi}(X_{H})$ into the Banach space $H_{p,\phi}(U_{i})$ of bounded 
$l_{p,\phi}(X_{H})$-valued holomorphic functions
$f$ on $U_{i}$ with norm
$$
|f|_{p,\phi}^{U_{i}}:=\sup_{z\in U_{i}}||f(z,\cdot)||_{p,\phi}
$$
(here $||\cdot||_{p,\phi}$ is the norm on $l_{p,\phi}(X_{H})$, see
(\ref{p1})).\\
(b) For every $z\in U_{i}\cap U_{j}$
$$
T_{i}(v)(z,\tau_{Q(H)}(\widetilde c_{ij})(g))=T_{j}(v)(z,g)\ ,\ \ \
v\in l_{p,\phi}(X_{H})\ .
$$
\end{Lm}
{\bf Proof.} (a) Let us first check the above statement for $p=1$ and
$p=\infty$. In fact, for $p=1$ from (\ref{oper}) and (\ref{bou1}) we have
$$
\begin{array}{c}
\displaystyle
||T_{i}(v)(z,\cdot)||_{1,\phi}:=\sum_{g\in X_{H}}
\left|\sum_{h\in X_{H}}v(h)S_{i_{0}i}(x,h,z,g)\right|\phi(g)\leq\\
\\
\displaystyle
\sum_{g\in X_{H}}
\sum_{h\in X_{H}}|v(h)\phi(h)|\cdot |S_{i_{0}i}(x,h,z,g)|
\frac{\phi(g)}{\phi(h)}
\leq\\
\\
\displaystyle
\left(\sum_{h\in X_{H}}|v(h)\phi(h)|\right)\cdot\sup_{h\in X_{H}}
\left(\sum_{g\in X_{H}}|S_{i_{0}i}(x,h,z,g)|\frac{\phi(g)}{\phi(h)}\right)\leq
\\
\\
\displaystyle
||v||_{1,\phi}\cdot |S_{i_{0}i}(x,\cdot,z,\cdot)|_{B}\leq C\cdot 
||v||_{1,\phi}\ .
\end{array}
$$
From here we have (see (\ref{se1}))
$$
|T_{i}(v)|_{1,\phi}^{U_{i}}\leq C\cdot ||v||_{1,\phi}\ .
$$
Thus $T_{i}(v)\in H_{1,\phi}(U_{i})$, and 
$T_{i}: l_{1,\phi}(X_{H})\to H_{1,\phi}(U_{i})$ is a linear bounded operator.
Similarly, for $p=\infty$ we have
$$
\begin{array}{c}
\displaystyle
|T_{i}(v)|_{\infty,\phi}^{U_{i}}:=\sup_{z\in U_{i}}
\left(\sup_{g\in X_{H}}\left|\sum_{h\in X_{H}}v(h)
S_{i_{0}i}(x,h,z,g)\right|\right)
\leq\\
\\
\displaystyle
||v||_{\infty,\phi}\cdot\sup_{z\in U_{i}}
\left(\sup_{g\in X_{H}}\left(\sum_{h\in X_{H}}
|S_{i_{0}i}(x,h,z,g)|\right)\right)
\leq 
C\cdot ||v||_{\infty,\phi}\ .
\end{array}
$$
So, $T_{i}: l_{\infty,\phi}(X_{H})\to H_{\infty,\phi}(U_{i})$ is well 
defined and continuous.

Let us prove the similar statement for $1<p<\infty$. 
Consider the evaluation $T_{z,i}(v)$ of $T_{i}(v)$ at $z\in U_{i}$\ :
$$
[T_{z,i}(v)](g):=T_{i}(v)(z,g)=\sum_{h\in X_{H}}v(h)S_{i_{0}i}(x,h,z,g)\ .
$$
From the above arguments it follows that $T_{z,i}$ is
a linear continuous map of $l_{1,\phi}(X_{H})$ to $l_{1,\phi}(X_{H})$ and of
$l_{\infty,\phi}(X_{H})$ to $l_{\infty,\phi}(X_{H})$, and in both these 
cases its
norm is bounded by $C$. Now, by the M. Riesz interpolation theorem
(see, e.g., [R]), $T_{z,i}$ maps also each $l_{p,\phi}(X_{H})$ to 
$l_{p,\phi}(X_{H})$
and its norm there is bounded by $C$, as well. Taken the supremum of norms
of $T_{z,i}$ over $z\in U_{i}$ we obtain that $T_{i}$ is a linear continuous 
operator from $l_{p,\phi}(X_{H})$ into $H_{p,\phi}(U_{i})$ for any $p$.\\
(b) Using (\ref{intersect}) with $i=i_{0}$, $j=i$ and $l=i_{0}$,
$m=j$ we get
$$
S_{i_{0}i}(x,h,z,\tau_{Q(H)}(\widetilde c_{ij})(g))=S_{i_{0}j}(x,h,z,g)\ .
$$
This and (\ref{oper}) produce the result.\ \ \ \ \ $\Box$

From Lemma \ref{bounded}  we obtain that
for every $v\in l_{p,\phi}(X_{H})$ the family $\{T_{i}(v)\}_{i\in I}$
represents a section $T(v)$ from $B_{p,\phi}(X_{H})$, see (\ref{se1}).
Moreover, the correspondence $v\mapsto T(v)$
determines a linear bounded operator 
$T:l_{p,\phi}(X_{H})\to B_{p,\phi}(X_{H})$. Note also that
(\ref{rest1}) implies that $T_{i_{0}}(v)(x,\cdot)=v$.
From this identifying $l_{p,\phi}(X_{H})$ with the fibre of 
$E_{p,\phi}(X_{H})$ over $x$ (by means of the coordinates on $U_{i_{0}}$) we 
obtain that $T$ maps the sections of $E_{p,\phi}(X_{H})|_{\{x\}}$ 
to $B_{p,\phi}(X_{H})$ and $T(v)(x)=v$.
To complete the proof of the proposition it remains to identify 
$B_{p,\phi}(X_{H})$ with ${\cal H}_{p,\psi}(M')$ and the space of
sections of $E_{p,\phi}(X_{H})|_{\{x\}}$ with $l_{p,\psi,x}(M')$.
\ \ \ \ \ $\Box$\\
{\bf 3.2.} We proceed with
the proof of the theorem. According to condition (\ref{e1}) there exists a
connected Stein neighbourhood $\widetilde M$ of $\overline{M}$ such that
$\widetilde M\subset\subset N$. Let
$\widetilde M'=r^{-1}(\widetilde M)$ be
the corresponding covering of $\widetilde M$. Let us consider the space
${\cal H}_{p,\psi}(\widetilde M')$. By 
$E_{0}(\widetilde M):=\widetilde M\times {\cal H}_{p,\psi}(\widetilde M')$ we 
denote the trivial holomorphic Banach vector bundle on $\widetilde M$ with 
fibre ${\cal H}_{p,\psi}(\widetilde M')$. Also, by $E_{p,\phi}(\widetilde M)$
we denote the bundle $E_{p,\phi}(X_{H})$ 
(see Example \ref{constbun} (b)). (Here $\phi$ is defined as in 
Proposition \ref{pr1}.)

For $z\in \widetilde M$ let $R_{z}$ be the restriction map of 
functions from ${\cal H}_{p,\psi}(\widetilde M')$ to the fibre $r^{-1}(z)$. 
If we identify ${\cal H}_{p,\psi}(\widetilde M')$ with
the Banach space $B_{p,\phi}(\widetilde M)$ of bounded holomorphic sections 
of the bundle $E_{p,\phi}(\widetilde M)$ by the direct image map 
(cf. Proposition \ref{sec1}), then $R_{z}$ will be the evaluation map of
sections from $B_{p,\phi}(\widetilde M)$ at $z\in\widetilde M$. In
particular, one can define a homomorphism of bundles
$R: E_{0}(\widetilde M)\to E_{p,\phi}(\widetilde M)$ which maps
$z\times v\in E_{0}(\widetilde M)$ to the vector $R_{z}(v)$ in the fibre
over $z$ of the bundle $E_{p,\phi}(\widetilde M)$ (see Definition
\ref{de2}).
Now, by Proposition \ref{pr1} with $M$ replaced by $\widetilde M$, 
every $R_{z}$ is surjective and, moreover,
there exists a linear continuous map $C_{z}$ of the fibre 
$E_{p,\phi, z}(\widetilde M)$ of $E_{p,\phi}(\widetilde M)$
over $z$ to the fibre $E_{0,z}(\widetilde M)$ of $E_{0}(\widetilde M)$ over 
$z$ such that $R_{z}\circ C_{z}=id$. Finally, by
$Ker\ \!R:=\cup_{z\in\widetilde M}Ker\ \!R_{z}\ (\subset E_{0}(\widetilde M))$
we denote the kernel of $R$.

Let $U\subset\widetilde M$ be an open set biholomorphic to
a Euclidean ball. Using some holomorphic
trivializations of $E_{0}(\widetilde M)|_{U}$ and 
$E_{p,\phi}(\widetilde M)|_{U}$ we may assume that
$E_{0}(\widetilde M)|_{U}=U\times B_{p,\phi}(\widetilde M)$ and 
$E_{p,\phi}(\widetilde M)|_{U}=U\times l$, where $l=l_{p,\phi}(X_{H})$ is 
the fibre of $E_{p,\phi}(\widetilde M)$.
In these trivializations we have $R(z\times v):=z\times R_{z}(v)$,
$z\times v\in U\times B_{p,\phi}(\widetilde M)$, where
$R_{z}(v):=v(z)$, and 
$C_{z}(z\times f):=z\times C_{z}(f)$, $z\times f\in U\times l$, where 
$C_{z}:l\to B_{p,\phi}(\widetilde M)$ is a linear continuous map. 

Next, take $x\in U$. Then in the above trivializations
$R_{z}\circ C_{x}:l\to l$, $z\in U$, is a 
family of linear continuous operators, holomorphic in $z$,
such that $R_{x}\circ C_{x}=id$.
Thus by the inverse function theorem (which in this case follows easily 
from the Taylor expansion of $R_{z}\circ C_{x}$ at $x$) we obtain
that
{\em there exists a neighbourhood $U_{x}$ of $x$ such that
$R_{z}\circ C_{x}$ is invertible for every $z\in U_{x}$.}

We set $P_{z}:=(R_{z}\circ C_{x})^{-1}$, and 
$\widetilde C_{z}:=C_{x}\circ P_{z}:l\to B_{p,\phi}(\widetilde M)$, 
$z\in U_{x}$. Then $P_{z}$ and $\widetilde C_{z}$ are holomorphic in 
$z\in U_{x}$, and by definition 
$$
R_{z}\circ\widetilde C_{z}=id\ .
$$ 
It is easy to see (by the inverse function theorem) that this 
identity implies that 
{\em there is a neighbourhood
$V_{x}\subset U_{x}$ of $x$ such that $Ker\ \!R|_{V_{x}}$ is biholomorphic
to $V_{x}\times Ker\ \!R_{x}$ and this biholomorphism is linear on
every $Ker\ \!R_{z}$ and maps this space onto $z\times Ker\ \!R_{x}$,
$z\in V_{x}$.}

Taking different $U$ and $x$ we then obtain from here that
$Ker\ \!R$ is a holomorphic Banach vector subbundle of 
$E_{0}(\widetilde M)$ and locally, for every $V_{x}$ as above, we have
an isomorphism of bundles
$E_{0}(\widetilde M)|_{V_{x}}\cong E_{p,\phi}(\widetilde M)|_{V_{x}}\oplus 
(Ker\ \! R)|_{V_{x}}$ given in the above trivializations by the formula
$$
z\times v\mapsto (z\times R_{z}(v), \
z\times (id-\widetilde C_{z}\circ R_{z})(v))\ .
$$

Let ${\cal U}=\{U_{i}\}_{i\in I}$ be an open cover
of $\widetilde M$ by Stein sets such that every quotient bundle 
$E_{p,\phi}(\widetilde M)|_{U_{i}}$ is complemented in 
$E_{0}(\widetilde M)|_{U_{i}}$. Let 
$C_{i}:E_{p,\phi}(\widetilde M)|_{U_{i}}\to 
E_{0}(\widetilde M)|_{U_{i}}$ be the complement homomorphism (i.e., 
$R|_{U_{i}}\circ C_{i}$ is the 
identity homomorphism of $E_{p,\phi}(\widetilde M)|_{U_{i}}$). We set 
$$
C_{ij}(z)=C_{i}(z)-C_{j}(z)\ ,\ \ \ z\in U_{i}\cap U_{j}\ .
$$
Then $\{C_{ij}\}$ is a holomorphic 1-cocycle with values in $Ker\ \!R$.
Since ${\cal U}$ is acyclic and $\widetilde M$ is Stein,
by the Bungart theorem [B] and the Leray lemma
one can find a family $\{H_{i}\}_{i\in I}$ of holomorphic sections of
$Ker\ \!R$ over $U_{i}$ such that
$$
H_{i}(z)-H_{j}(z)=C_{ij}(z)\ ,\ \ \ z\in U_{i}\cap U_{j}\ .
$$
In particular, setting $F|_{U_{i}}:=C_{i}-H_{i}$ we obtain that
$F:E_{p,\phi}(\widetilde M)\to E_{0}(\widetilde M)$ is a (holomorphic)
homomorphism of bundles such that $R\circ F=id$. 
Since by the definition $\overline{M}\subset\widetilde M$ is compact, 
$F|_{M}$ is bounded. Next, by $R_{M}$ we denote the restriction 
homomorphism of trivial bundles $E_{0}(\widetilde M)|_{M}\to E_{0}(M)$
defined by 
$$
z\times v\mapsto z\times v|_{M}\ ,\ \ \ z\in M\ ,\ 
v\in {\cal H}_{p,\psi}(\widetilde M')\ .
$$
Finally, we set 
$$
L_{z}:=R_{M}\circ F(z)\ ,\ \ \ z\in M\ .
$$
Then by definition, every $L_{z}$ is a linear continuous map of
$l_{p,\psi,z}(M')$ into ${\cal H}_{p,\psi}(M')$ (see section 1.1),
the family $\{L_{z}\}$ is holomorphic in $z\in M$, 
$R_{z}\circ L_{z}=id$, and $\sup_{z\in M}||L_{z}||<\infty$.

This completes the proof of the theorem.\ \ \ \ \ $\Box$
\begin{R}\label{unif}
{\rm Using some modification of the arguments of the above proof one can 
show that for a fixed $\psi$ there exists a family $\{L_{z}\}$ satisfying
the assumptions of Theorem \ref{te1} such that
$\sup_{z\in M}||L_{z}||\leq C<\infty$ where $C$ does not depend
on the class ${\cal H}_{p,\psi}(M')$ (cf. the construction in the proof of
Theorem 1.4. of [Br2]).}
\end{R}
%==============
\sect{\hspace*{-1em}. Proofs.}
{\bf Proof of Corollary \ref{co1}.} The result follows straightforwardly
from formula (\ref{le1}) where $f$ is a holomorphic function with 
values in ${\cal H}_{p,\psi}(M')$ and from the properties of the operators
$L_{z}$.\ \ \ \ \ $\Box$\\
{\bf Proof of Theorem \ref{te2}.}
Let $f\in {\cal H}(r^{-1}(C))$ be a function satisfying the assumptions of the
theorem. Consider the function
$$
h(z):=L_{z}(f|_{r^{-1}(z)})\ ,\ \ \ z\in C\ .
$$
By Definition \ref{cont1}, Proposition \ref{sec1} and by the properties of 
$\{L_{z}\}$ we have
that $h$ is a ${\cal H}_{p,\psi}(M')$-valued continuous function on $C$ 
holomorphic in $D$. (It can be written as the scalar function of the
variables $(z,w)\in C\times M'$.) Therefore it suffices to prove
an approximation theorem for such Banach-valued functions. 
Namely, it suffices to show that 
{\em every such $h$ can be uniformly approximated on
$C$ by ${\cal H}_{p,\psi}(M')$-valued holomorphic functions
defined in a neighbourhood $\Omega$ of $C$}. 

Further, if $\{h_{i}(z,w)\ :\ z\in\Omega,\ w\in M'\}_{i\geq 1}$ is such an
approximation sequence for $h$, then 
$\{f_{i}(y)\ :\ y\in r^{-1}(\Omega)\}_{i\geq 1}$ with
$f_{i}(y):=h_{i}(r(y),y)$ is the approximation sequence
for $f$ satisfying the required statement of the theorem.

The proof of the above approximation theorem for functions $h$ repeats
word-for-word the proof of Theorem 3.5.1 in [HL] where in all
integral formulas we replace the scalar functions by Banach-valued ones.
We leave the details to the readers.\ \ \ \ \ $\Box$\\
{\bf Proof of Theorem \ref{te3}.} Let us consider a function $f$ satisfying
either (1) or (2). Define the function
$$
h(z):=L_{z}(f|_{r^{-1}(z)})\ ,\ \ \ z\in X\ ({\rm or}\ \ z\in\overline{X})\ .
$$
Then according to Definition \ref{cont1}, $h$ is a holomorphic function on 
the submanifold $X\subset M$ with values in 
${\cal H}_{p,\psi}(M')$ (and in case (2) it is also continuous on 
$\overline{X}$).
Thus it suffices to prove the extension theorem for Banach-valued
holomorphic functions on $X$ (extending them to $M$). Evaluating
the extended Banach-valued functions at the points $(r(y),y)$, $y\in M'$ 
(cf. the proof of Theorem \ref{te2} above), we get the required result. 

The scalar case of the required extension theorem is proved in 
Theorem 4.11.1 of [HL]. The Banach-valued case repeats literally
the arguments of the proof of Theorem 4.11.1 
where in all integral formulas we replace  scalar functions by
Banach-valued ones. Also, instead of the
classical Cartan B theorem for Stein manifolds, we use in the 
proof its Banach-valued generalization due to Bungart.\ \ \ \ \ $\Box$
\begin{R}\label{boun1}
{\rm One can show (cf. the remark after the proof of Theorem 4.11.1 of [HL])
that under the hypothesis (1) of Theorem \ref{te3} 
there exists a linear continuous operator 
$E: {\cal H}_{p,\psi}(r^{-1}(X))\to {\cal H}_{p,\psi}(M')$ such that
$Ef=f$ on $r^{-1}(X)$ for all $f\in {\cal H}_{p,\psi}(r^{-1}(X))$.
Moreover, the norm $||E||$ of $E$ is bounded by a constant
depending only on $\psi:M'\to\Re_{+}$ and $X$.
In particular, if $\psi\equiv 1$, then $||E||$ depends only on $X$
(and not on the covering $M'$). However, it is not clear
whether the extension in Theorem \ref{te3} (2) can be made by a linear
continuous operator, as well.}
\end{R}
%===============
\sect{\hspace*{-1em}. Appendix.}
We present a multi-dimensional analog of Cauchy-Green
formulas on coverings of domains $M\subset\Co^{n}$ satisfying condition
(\ref{e1}). These formulas are
obtained from similar ones on the domains $M$ (see, e.g.,
[H]) by the application of Theorem \ref{te1}. Below we use the notation
of the Leray integral formula (\ref{le1}) (see section 1.2).

Let $M\subset\Co^{n}$ be a domain with a rectifiable boundary. Assume that
$M$ satisfies condition (\ref{e1}), that is, 
$M\subset\subset\widetilde M\subset N\subset\Co^{n}$,
$\pi_{1}(N)=\pi_{1}(M)$ and $\widetilde M$ is Stein. 
Let $\eta=\eta(\xi,z)=(\eta_{1},\dots,\eta_{n})$ be a smooth 
$C^{n}$-valued function of the variable $\xi\in\overline{M}$ such that
$<\eta(\xi,z),\ \xi-z>=1$ for $\xi\in\partial M$. 
Consider a covering $M'$ of $M$. Then $M'$ is an open subset of the
covering $r: N'\to N$ such that $\pi_{1}(N')=\pi_{1}(M')$.
Let $\psi:N'\to\Re_{+}$ be such that $\log\psi$ is
uniformly continuous with respect to the path metric obtained by a
Riemannian metric pulled back from $N$. By $E_{p,\phi}(\widetilde M)$
we denote the bundle $E_{p,\phi}(X_{H})$ on $\widetilde M$
(see Example \ref{constbun} (b)). Here $\phi$ is defined
as in Proposition \ref{pr1}. Also, by $L=\{L_{z}\}$ we denote the 
family of operators given by Theorem \ref{te1} (with the same $p$ and
$\psi$) where $z$ varies in a small neighbourhood $\Omega$ of 
$\overline{M}$ (in this case $L_{z}$ maps 
$l_{p,\psi,z}(r^{-1}(\Omega))$ into 
${\cal H}_{p,\psi}(r^{-1}(\Omega)$)).

Suppose that a function $f$ defined on $\overline{M'}$ is such that its
direct image $r_{*}(f)$ with respect to $r$ is a continuous section of
$E_{p,\phi}(\widetilde M)|_{\overline{M}}$, and $\overline\partial r_{*}(f)$
is a continuous section of $E_{p,\phi}(\widetilde M)|_{\overline{M}}$.
We set
$$
(Lf)(z):=L_{z}(f|_{r^{-1}(z)})\ ,\ \ \ z\in\overline{M}\ .
$$
Then from the construction of the family $L$ and from the definition of
$f$ it follows that $Lf$ is a continuous  
${\cal H}_{p,\psi}(r^{-1}(\Omega))$-valued function on $\overline{M}$ 
such that $\overline\partial Lf$ is a continuous  
${\cal H}_{p,\psi}(r^{-1}(\Omega))$-valued $(0,1)$-form on $\overline{M}$.

Now, for the function $f$ we have the 
following integral representations:
\begin{equation}\label{last1}
f(z)=\widetilde K^{s}f(z)+\widetilde H^{s}(\overline\partial f)(z)\ ,
\ \ \ z\in M'\ ,
\end{equation}
where, for $y:=r(z)$,
$$
\widetilde K^{0}f(z)=\frac{(n-1)!}{(2\pi i)^{n}}\int_{\xi\in\partial M}
L_{\xi}(f|_{r^{-1}(\xi)})(z)\ \!\omega'(\eta(\xi,y))\wedge\omega(\xi)\ ,
$$
$$
\begin{array}{l}
\displaystyle
\widetilde H^{0}(\overline\partial f)(z)=\\
\\
\displaystyle
\frac{(n-1)!}{2\pi i)^{n}}\left(
\int_{(\xi,\lambda_{0})\in\partial M\times [0,1]}
(\overline\partial L_{\xi}(f|_{r^{-1}(\xi)}))(z)
\wedge\omega'\left((1-\lambda_{0})
\frac{\overline\xi-\overline y}{|\xi-y|^{2}}+\lambda_{0}\eta(\xi,y)\right)
\wedge\omega(\xi)\right.\\
\\
\displaystyle
+\left.\int_{\xi\in M}(\overline\partial L_{\xi}(f|_{r^{-1}(\xi)}))(z)
\wedge\omega'\left(
\frac{\overline\xi-\overline y}{|\xi-y|^{2}}\right)\wedge\omega(\xi)\right)\ .
\end{array}
$$
And for $s>0$,
$$
\begin{array}{l}
\displaystyle
\widetilde K^{s}f(z)=\\
\\
\displaystyle
\frac{(n+s-1)!}{(s-1)!(2\pi i)^{n}}\int_{\xi\in M}
L_{\xi}(f|_{r^{-1}(\xi)})(z)(1-<\eta(\xi,y),\ \xi-y>)^{s-1}\cdot
\omega(\eta(\xi,y))\wedge\omega(\xi),
\end{array}
$$
$$
\begin{array}{l}
\displaystyle
\widetilde H^{s}(\overline\partial f)(z)=\\
\\
\displaystyle
\frac{(n+s-1)!}{(s-1)!(2\pi i)^{n}}
\int_{(\xi,\lambda_{0})\in M\times [0,1]}
(\overline\partial L_{\xi}(f|_{r^{-1}(\xi)}))(z)\times
\left(\lambda_{0}(1-<\eta(\xi,y),\ \xi-y>\right)^{s-1}\times\\
\\
\displaystyle
\omega'\left((1-\lambda_{0})
\frac{\overline\xi-\overline z}{|\xi-y|^{2}}+\lambda_{0}\eta(\xi,y)\right)
\wedge\omega(\xi)\ .
\end{array}
$$
%====================

\end{document}